\documentclass[twoside,12pt]{amsart} 
\usepackage[dvips]{graphicx}
\DeclareGraphicsExtensions{ps,eps}
\makeatletter
\def\serieslogo@{}
\makeatother
\makeatletter 
\def\@setcopyright{}
\makeatother

\newtheorem{Theorem}{Theorem}[section]

\newtheorem{Definition}[Theorem]{Definition}
\newtheorem{Proposition}[Theorem]{Proposition}

\theoremstyle{definition}

\theoremstyle{remark}
\newtheorem{rem}{Remark}[section]

\begin{document}

\title
{
The initial value problem for the Euler equations of incompressible fluids
viewed as a concave maximization problem 
}

\author{Yann Brenier}
\address{
CNRS, Centre de math\'ematiques Laurent Schwartz,
Ecole Polytechnique, FR-91128 Palaiseau, France.
}
\curraddr{}
\email{yann.brenier@polytechnique.edu}

\maketitle
\markboth{}{}

\section*{Abstract}
We consider the Euler equations of incompressible fluids \cite{AK,Li}
and attempt to solve the initial value problem with the help of a concave maximization problem.
We show that this problem, which shares a similar structure with the optimal transport problem
with quadratic cost, in its "Benamou-Brenier" formulation \cite{BB,AGS,Vi},
always admits a relaxed solution that can be interpreted
in terms of $sub-solution$ of the Euler equations in the sense of convex integration theory
\cite{DS}.
Moreover, any smooth solution of the Euler equations can be recovered  
from this maximization problem, at least for short times.
\subsection*{Keywords}
Fluid mechanics, partial differential equations, calculus of variation, optimal transport,
Euler equations
\section*{Introduction}
Let us fix a time interval $[0,T]$ 
and denote by $D$ the periodic box $D=\mathbb{R}^d/\mathbb{Z}^d$.
The Euler model 
\cite{AK,Li} 
of an incompressible fluid of unit mass density,
moving in $D$ during the time interval $[0,T]$, without
external forces, assumes the existence of a square-integrable, divergence-free vector field $V$, over $[0,T]\times D$,
such that (written in coordinates,
with the usual implicit summation
over repeated lower and upper indices)
$
\partial_t V^i+\partial_j(V^i V^j)
$
is a gradient field. In weak form, this means 
\begin{equation}
\label{weak0}
\int_{[0,T]\times D}
\partial_i\varphi(t,x)V^i(t,x)dxdt=0,
\end{equation}
for all smooth function $\varphi$ over $[0,T]\times D$,
which encodes that $V$ is divergence-free and 
\begin{equation}
\label{weak}
\int_{[0,T]\times D}
\left(\partial_jA_iV^i V^j+\partial_t A_iV^i\right)(t,x)dxdt
+\int_D P_0^i(x)A_i(0,x)dx=0,
\end{equation}
for all smooth divergence-free vector fields $A$
on $[0,T]\times D$,
vanishing at $t=T$, 
which includes (weakly) the initial condition that $V$ is $P_0$ at time $t=0$,
$P_0$ being a given $L^2$ divergence-free vector field on $D$.

Our goal is to solve, by a concave maximization method, the initial value problem for the Euler model with as initial condition a 
fixed divergence-free vector field $P_0$, square integrable over $D$ and of zero spatial mean.
The idea is very simple: we try to find a divergence-free vector field $V$, 
weak solution to the Euler equation with initial condition $P_0$, of minimal kinetic energy. This leads to the saddle-point problem

\begin{equation}
\label{ivp0}
\mathcal{I}[P_0]
=\inf_V \sup_{A,\varphi}
\;\int_{[0,T]\times D}
\frac{1}{2}{|V|}^2
+\partial_jA_iV^i V^j+(\partial_t A_i+\partial_i\varphi)V^i
+\int_D P_0^iA_i(0,\cdot)
\end{equation}
$$
$$
over all  $L^2$ vector fields $V$
on $[0,T]\times D$,
all smooth divergence-free vector fields $A$ vanishing 
at $t=T$, and all smooth
real functions $\varphi$. We may interpret $(A,\varphi)$ as Lagrange multipliers 
for the constraint that $V$ is a weak solution to the Euler equations with initial condition $P_0$, 
in the sense of (\ref{weak0},\ref{weak}).
Investigating problem (\ref{ivp0}) looks silly since the Euler equation to be solved is already
included as a constraint! Furthermore, for 
smooth solutions of the Euler equation on the periodic box,
the kinetic energy is constant in time and, therefore, depends only on the data $P_0$,
so that...there seems to be nothing to minimize! However, for a fixed initial condition,
weak solutions are not unique and the conservation of energy is generally not true
as well known since the celebrated results of Scheffer, Shnirelman,
De Lellis and Sz\'ekelyhidi \cite{Sc,Sh,DS}. Therefore, since, in addition, weak solutions always exist, following Wiedemann \cite{Wi}, the minimization problem is definitely not meaningless.
\\
\\
In this paper, we mostly investigate the dual problem obtained by 
exchanging the infimum and
the supremum in (\ref{ivp0}), leading to a concave maximization problem which will be
shown to be solvable in section \ref{existence0}, after a suitable reformulation of the concept
of weak solutions established in section \ref{reformulation}.
The resulting maximization problem roughly reads
$$
\sup_{(E,B)}
\;\;-
\int_{[0,T]\times D}
E\cdot(\mathbb{I}_d+2B)^{-1}\cdot E+2P_0\cdot E
$$
where $\mathbb{I}_d$ denotes the $d\times d$ identity matrix, $E$ and $B$ are respectively valued  in $\mathbb{R}^d$ and in the space of $d\times d$ symmetric
matrices, and subject to
$$
\partial_t B=LE,\;\;\;B(t=T,\cdot)=0,
$$
$L$ being a suitable first-order constant coefficient (pseudo-)differential operator on $D$,
namely (in coordinates)
$$
L^k_{ij}E_k=
\frac{1}{2}(\partial_j E_i+\partial_i E_j)
+\partial_i\partial_j(-\bigtriangleup)^{-1}\partial^k E_k,\;{\rm{where}}\;\partial^k=\delta^{kj}\partial_j\;
{\rm{and}}\;\bigtriangleup=\delta^{ij}\partial_i\partial_j.
$$
Surprisingly enough, this problem looks very similar to the Monge optimal mass transport problem with quadratic cost in its "Benamou-Brenier"
formulation \cite{BB,AGS,Vi}, which would read
$$
\inf_{\rho,Q} \int_{[0,T]\times D}
Q\cdot\rho^{-1}\cdot Q,
$$
where $\rho$ and $Q$ are respectively valued in $\mathbb{R}_+$ and $\mathbb{R}^d$ and subject to $\partial_t\rho+\partial_i Q^i=0$, while $\rho$ is prescribed at $t=0$ and $t=T$. Presumably, the maximization problem can be treated by the same numerical method as the one used in \cite{BB}.
\\
\\
Next, in section \ref{consistency}, we will
establish the consistency result that any local smooth solution of the Euler equations can be recovered,
from the maximization problem, for short enough times $T$. 
\\
\\
Finally, in section \ref{duality0}, we will make a connection
between the maximization problem and the theory of sub-solutions to the Euler equations
which has recently attracted a lot of interest after the celebrated work of De Lellis and
Sz\'ekelyhidi \cite{DS} in the framework of Convex Integration theory.
\subsection*{Acknowledgment}
This work has been partly supported by the grant ISOTACE
ANR-12-MONU-0013 (2012-2016). The author thanks Nassif Ghoussoub
for exciting discussions (held at the Erwin Schr\"odinger Institute in the
summer of 2016) about his new theory of "optimal ballistic transport"
\cite{Gh}, which were certainly influential for the present paper.
He is also grateful to the Schr\"odinger Institute and to the
CNRS-INRIA MOKAPLAN team for their help during the preparation of this work.
\section{Reformulation of the weak formulation of the Euler equations}
\label{reformulation}
We first revisit definition (\ref{weak0},\ref{weak}) of weak
solutions to the Euler equations by substituting for $(A,\varphi)$
the smooth fields $(E,B)$
respectively valued
in $\mathbb{R}^d$ and in the set of 
symmetric $d\times d$ matrices, 
defined (in coordinates) on $[0,T]\times D$ by
$$
E_i=\partial_t A_i+\partial_i\varphi\;,
\;\;\;B_{ij}=
\frac{1}{2}(\partial_jA_i+\partial_i A_j).
$$
We first observe that $B(T,\cdot)=0$, while
$(E,B)$ satisfy the compatibility condition:
$$
\partial_t B_{ij}=
\frac{1}{2}(
\partial_jE_i+\partial_i E_j)
-\partial_i\partial_j\varphi.
$$
Since $A$ is divergence free, $B$ is trace-less and, therefore,
$$
\bigtriangleup \varphi=\partial^k E_k,\;\;\;{\rm{with\;notations:}}\;\partial^k=\delta^{jk}\partial_j,\;\;\;
\bigtriangleup \varphi=\partial^k \partial_k.
$$
Thus, in short,
$(E,B)$ are just subject to 
constraint
\begin{equation}
\label{constraint-B-E}
\partial_t B=LE,\;\;\;B(T,\cdot)=0.
\end{equation}
where $L$ is the constant (pseudo-)differential operator of order 1 defined by
\begin{equation}
\label{L}
L^k_{ij}E_k=
\frac{1}{2}(\partial_j E_i+\partial_i E_j)
+\partial_i\partial_j(-\bigtriangleup)^{-1}\partial^k E_k,\;{\rm{where}}\;\partial^k=\delta^{kj}\partial_j\;
{\rm{and}}\;\bigtriangleup=\delta^{ij}\partial_i\partial_j.
\end{equation}
Constraint (\ref{constraint-B-E}) can be, as well, written in weak form:
\begin{equation}
\label{constraint-B-E+}
\int_{(t,x)\in[0,T]\times D}
B_{ij}(t,x)\zeta'(t)\psi^{ij}(x)
-E_k(t,x)\zeta(t)
\left(\partial_j\psi^{kj}-\partial^k(-\bigtriangleup)^{-1}\partial_i\partial_j\psi^{ij}\right)(x)=0
\end{equation}
for all smooth functions $\zeta$ on $[0,T]$ vanishing at $t=0$ and all smooth functions $\psi$ on $D$ 
valued in the set of $d\times d$ symmetric matrices.
Notice that expression (\ref{constraint-B-E+}) still makes sense
as $(E,B)$ are just bounded Borel measures over $[0,T]\times D$.
Later on, we will use the class, denoted by $\mathcal{EB}$,
of such measures subject to (\ref{constraint-B-E+}).
\\
\\
Next, in order to reformulate (\ref{weak0},\ref{weak})
entirely in terms of $(E,B)$ instead of $(A,\varphi)$
we have to express the time-boundary term
$
\int_D P_0\cdot A(0)
$ 
in terms of $(B,E)$. We have
$$
A(0,x)=-\int_0^T \partial_t A(t,x)dt,
$$ 
since $A(T,\cdot)=0$. Thus
$$
\int_D P_0\cdot A(0)=-\int_{[0,T]\times D} P_0\cdot\partial_t A=-\int_{[0,T]\times D} P_0\cdot E
$$ 
(by definition of $E$, using that $P_0$ is divergence-free 
with zero spatial mean).
Finally, using $(E,B)$ instead of $(A,\varphi)$, in definitions
(\ref{weak0},\ref{weak}), we conclude:
\begin{Proposition}
\label{weak-new}
Given a divergence-free zero-mean $L^2$-vector-field $P_0$ on $D$,
an $L^2$ vector-field  $V$ over $[0,T]\times D$ is a weak solution to the Euler equations with initial condition $P_0$, in the sense of 
(\ref{weak0},\ref{weak}), if and only if
\begin{equation}
\label{weakbis}
\int_{[0,T]\times D} V^i V^j B_{ij}+(V^i-P_0^i)E_i=0,
\end{equation}
for all smooth fields $(E,B)$
respectively valued
in $\mathbb{R}^d$ and in the set of 
symmetric $d\times d$, that are subject to constraint (\ref{constraint-B-E}),
namely
$$
B(T,\cdot)=0,\;\;\;\partial_t B_{ij}=
\frac{1}{2}(\partial_j E_i+\partial_i E_j)
+\partial_i\partial_j(-\bigtriangleup)^{-1}\partial^k E_k.
$$
In addition, in the new definition (\ref{weakbis}) of weak solutions, 
we may extend the range of trial fields $(E,B)$ to the class 
$\mathcal{EB}_{2,\infty}$ which is the intersection of $L^2\times L^\infty$
with the class $\mathcal{EB}$ of all bounded Borel measures  on $[0,T]\times D$
satisfying (\ref{constraint-B-E+}).
\end{Proposition}
The proof is straightforward: the first statement directly follows from the previous calculations,
while the second one follows from the rather obvious property that $\mathcal{EB}_{2,\infty}$ is just the $L^2\times L^\infty$ weak-* 
closure of the class $\mathcal{EB}_{smooth}$ of its own smooth elements. 
$[$Indeed, every $(E,B)$ in $\mathcal{EB}_{2,\infty}$ can be approximated, as closely as needed,
by $(\tilde E,\tilde B)\in\mathcal{EB}_{smooth}$ in two steps: we first,
mollify $E$ on $[0,T]\times D$ and get $\tilde E$; then, we just set $\tilde B(t,x)=-\int_t^T L\tilde E(s,x)ds$ so that 
(\ref{constraint-B-E}), and therefore (\ref{constraint-B-E+}), is still satisfied by $(\tilde E,\tilde B)$, while $\tilde B$
stays close to $B$.$]$

\section{A concave maximization problem}
\label{existence0}
Using the new definition (\ref{weakbis}) of weak solutions we found in the previous section, we may formulate
the original minimization problem (\ref{ivp0}) just as
\begin{equation}
\label{ivp1}
\mathcal{I}[P_0]
=\inf_V \sup_{(E,B)\in\mathcal{EB}_{2,\infty}}
\;\;\int_{[0,T]\times D} \frac{1}{2}V^i V^j (\delta_{ij}+2B_{ij})+(V^i-P_0^i)E_i,
\end{equation}

Since $\inf\sup\ge \sup\inf$, we get the lower bound 
$$
\mathcal{I}[P_0]\ge \mathcal{J}[P_0]
=\sup_{(E,B)\in\mathcal{EB}_{2,\infty}}\inf_V 
\;\;\int_{[0,T]\times D} \frac{1}{2}V^i V^j (\delta_{ij}+2B_{ij})+(V^i-P_0^i)E_i.
$$
In this "dual" problem, the infimum over all $L^2$ vector fields $V$ on
${[0,T]\times D\times \mathbb{R}^d}$ is very easy to deal with. Indeed, it is certainly equal to 
$-\infty$ unless 
\begin{equation}\label{vp}
\mathbb{I}_d+2B\ge 0,\;\;\;a.e.\;\;{\rm{in}}\; [0,T]\times D,
\end{equation}
in the sense of symmetric matrices. 
Next, we observe that, because $(E,B)$ belongs to $\mathcal{EB}_{2,\infty}$, $B$ must be trace-free.
Indeed, this follows from (\ref{constraint-B-E}), using definition (\ref{L}) of $L$. So, denoting by 
$\lambda(\alpha)\in\mathbb{R}$ the eigenvalues of $B$, 
for $\alpha\in\{1,\cdot\cdot\cdot,d\},$ we have (using (\ref{vp})):
$$
0\le 1+2\lambda(\alpha)\le \sum_{\beta=1}^{d}(1+2\lambda(\beta))=d,
$$
and, therefore,
\begin{equation}\label{vpbound0}
-\mathbb{I}_d \le 2B\le (d-1)\mathbb{I}_d,
\end{equation}
which provides an a priori $L^\infty$ bound for B (since $B$ is valued in the set of symmetric
matrices).
We also get
\begin{equation}\label{vpbound}
(\mathbb{I}_d+ 2B)^{-1}\ge d^{-1}\mathbb{I}_d.
\end{equation}
and, after minimization in $V$,
\begin{equation}\label{ivpdual0}
\mathcal{J}[P_0]
=\sup_{(E,B)\in\mathcal{EB}_{2,\infty}}
\;\;-K[E,B]-\int_{[0,T]\times D}
P_0\cdot E,
\end{equation}
where
\begin{equation}\label{K0}
K[E,B]=\frac{1}{2}\int_{[0,T]\times D}E\cdot(\mathbb{I}_d+2B)^{-1}\cdot E.
\end{equation}
This immediately implies $J[P_0]\ge 0$
(just by taking $E=B=0$).
\\
Here, we emphasize that (\ref{ivpdual0}) is a concave maximization problem in $(E,B)$.
Indeed, we can write, point-wise in $(t,x)$,
$$
E\cdot(\mathbb{I}_d+2B)^{-1}\cdot E
=\sup_{M,Z} \;\;2E_iZ^i-(\delta_{ij}+2B_{ij})M^{ij}
$$
where $M$ and $Z$ are respectively $d\times d$ symmetric
matrices and vectors in $\mathbb{R}^d$ subject to
\begin{equation}\label{matrineq}
Z\otimes Z\le M,
\end{equation}
in the sense of symmetric matrices. This allows us to give a more precise definition of $K$, namely
\begin{equation}\label{K}
K[E,B]=\sup_{M\ge Z\otimes Z} \;\;\frac{1}{2}\int_{[0,T]\times D}
2 E_iZ^i-(\delta_{ij}+2B_{ij})M^{ij}
\in [-\infty,0],
\end{equation}
where the infimum is performed over all pairs $(Z,M)$ 
of continuous functions on
$[0,T]\times D$, respectively valued
in $\mathbb{R}^d$ and in the set of 
symmetric matrices $M$.
Notice that definition (\ref{K}) makes sense already as $(E,B)$ is just a pair of bounded Borel measures on $[0,T]\times D$, and a fortiori as $(E,B)$ belong to $L^2\times L^\infty$. In both cases,
$K$ is a lower semi-continuous convex function of $(E,B)$ valued in $[0,+\infty]$.
Next, because of the lower bound (\ref{vpbound}), we get an $L^2$ a priori bound for $E$. Indeed,
by definition (\ref{K0}) of $K$, we have:
$$
\frac{1}{2d}\int_{[0,T]\times D}|E|^2\le K[E,B].
$$
So, for any $\varepsilon-$maximizer $(E,B)\in \mathcal{EB}_{2,\infty}$ of (\ref{ivpdual0}), we get
$$
0\le \mathcal{J}[P_0]
\le \varepsilon-K[E,B]-\int_{[0,T]\times D}P_0\cdot E
\le \varepsilon-\frac{1}{4d}\int_{[0,T]\times D}(2|E|^2+4d P_0\cdot E)
$$
$$
= \varepsilon-\frac{1}{4d}\int_{[0,T]\times D}(|E|^2+|E+2dP_0|^2-4d^2|P_0|^2)
\le \varepsilon-\frac{1}{4d}\int_{[0,T]\times D}(|E|^2-4d^2|P_0|^2)
$$
which provides the a priori bound, for every $\varepsilon-$maximizer $(E,B)\in \mathcal{EB}_{2,\infty}$ of (\ref{ivpdual0}), 
\begin{equation}\label{bound-E}
\int_{[0,T]\times D}|E|^2\le 4d\varepsilon+4d^2 T\int_{D}|P_0|^2.
\end{equation}
We also deduce
$$
0\le \mathcal{J}[P_0]\le Td\int_{[0,T]\times D}|P_0|^2.
$$
$$
[{\rm{Indeed}}\;\;\;
0\le \mathcal{J}[P_0]\le\varepsilon-\frac{1}{4d}\int_{[0,T]\times D}(|E|^2-4d^2|P_0|^2)
\le\varepsilon+Td\int_{D}|P_0|^2.]
$$
By definition (\ref{K}), $K$ is lower semi-continuous with respect to the weak-* topology of 
$L^\infty \times L^2$, while
$$
E\rightarrow \int_{[0,T]\times D}P_0\cdot E,
$$
is continuous (since $P_0$ is given in $L^2$) and $\mathcal{EB}_{2,\infty}$ is weak-*
closed in $L^\infty \times L^2$.
Thus, we conclude that the maximization problem
(\ref{ivpdual0}) 
always has at least an optimal solution $(E,B)$ in class $\mathcal{EB}_{2,\infty}$, since its 
$\varepsilon-$maximizers stay confined in a fixed ball (and therefore a 
weak-* compact subset) 
of $L^\infty \times L^2$, as $\varepsilon$ goes to zero. 
\\
\\
In addition, let us observe that, for any $(E,B)$ in $\mathcal{EB}_{2,\infty}$,
$$
\int_{D} \left(B_{ij}(t_1,\cdot)-B_{ij}(t_0,\cdot)\right)\psi^{ij}\le
\sqrt{t_1-t_0\;}\;
||E||_{L^2([0,T]\times D)}||L^*\psi||_{L^\infty(D)}
$$
for all $t_0\le t_1$ in $[0,T]$ and or all smooth functions $\psi$ on $D$, valued in the set of $d\times d$ symmetric matrices, where
$$
(L^*)^k_{ij}\psi^{ij}
=\partial_j\psi^{kj}-\partial^k(-\bigtriangleup)^{-1}\partial_i\partial_j\psi^{ij}
$$
(because of (\ref{constraint-B-E}) and, more precisely,
(\ref{constraint-B-E+})). Thus
any maximizer $(E,B)$ must satisfy
$$
\int_{D} \left(B_{ij}(t_1,\cdot)-B_{ij}(t_0,\cdot)\right)\psi^{ij}
\le 2d\sqrt{T}
\sqrt{t_1-t_0\;}\;||P_0||_{L^2(D)}||L^*\psi||_{L^\infty(D)},
$$
because of (\ref{bound-E}).
This shows that $B$ belongs to $C^{1/2}([0,T],C^1(D)'_{w*})$.
\\
\\
So, we have finally obtained:
\begin{Theorem} 
\label{existence}
Let $P_0\in L^2(D,\mathbb{R}^d)$ be 
a divergence-free vector field 
of zero spatial mean.
Let $\mathcal{EB}_{\infty,2}$ be the class of all $L^2\times L^\infty$ fields $(E,B)$ on
$[0,T]\times D$,
respectively valued
in $\mathbb{R}^d$ and in the set of 
symmetric $d\times d$ matrices, that are subject to constraint 
$$
B(T,\cdot)=0,\;\;\;\partial_t B_{ij}=
\frac{1}{2}(\partial_j E_i+\partial_i E_j)
+\partial_i\partial_j(-\bigtriangleup)^{-1}\partial^k E_k
$$
(or, more precisely, (\ref{constraint-B-E+})).
Then the maximization problem 
$$
\mathcal{J}[P_0]
=\sup_{(E,B)\in\mathcal{EB}_{2,\infty}}
\;\;-
\int_{[0,T]\times D}\left(
\frac{1}{2}E\cdot(\mathbb{I}_d+2B)^{-1}\cdot E+P_0\cdot E\right)
$$
(or, more precisely, (\ref{ivpdual0},\ref{K})), always admits a solution $(E,B)$.
In addition,
$$
{\rm{spect}}(B)\subset[-1/2,(d-1)/2],\;\;
\int_{[0,T]\times D}|E|^2\le 4d^2 T\int_{D}|P_0|^2,\;\;
0\le\mathcal{J}[P_0]\le Td\int_{D}|P_0|^2,
$$
and $B$ belongs to $C^{1/2}([0,T],C^1(D)'_{w*})$.
\end{Theorem}
\rem
Notice that all the a priori bounds found for problem (\ref{ivpdual0}) 
are still valid when addressing the $relaxed$ maximization
problem
\begin{equation}\label{ivpdual}
\sup_{(E,B)\in\mathcal{EB}}
\;\;-K[E,B]-\int_{[0,T]\times D}
P_0\cdot E,
\end{equation}
where we only require $(E,B)$ to be plain bounded Borel measures and no longer in
the space $L^2\times L^\infty$, while $K$ should be understood
as in (\ref{K}) and $P_0$ is restricted to be in the class of continuous divergence-free vector fields,
with zero-mean, on $D$, so that
$$
E\rightarrow \int_{[0,T]\times D}P_0\cdot E,
$$
is still well defined and continuous.
Since the $\varepsilon-$maximizers of this relaxed problem 
must satisfy the $same$
a priori bounds (\ref{vpbound0},\ref{bound-E}) as the ones of  (\ref{ivpdual0}),
they $necessarily$ belong to the subspace $L^2\times L^\infty$, which implies that 
both problems (\ref{ivpdual0}) and (\ref{ivpdual}) 
admit the same optimal value $\mathcal{J}[P_0]$
and the same maximizers.
\\
In the last section of this paper, the relaxed problem (\ref{ivpdual})
will be reformulated and solved (at least in the case when $P_0$ is continuous)
by a duality method which will establish a link with the concept of
sub-solution used by De Lellis and Sz\'ekelyhidi in their approach of the Euler equations by
"Convex Integration" \cite{DS}.
\section{Recovery of smooth classical solutions to the Euler equations}
\label{consistency}
We want now to show that the optimization problem (\ref{ivpdual0}) addressed
in the previous section
is consistent with the classical theory of local smooth solutions of the initial value problem for the Euler equations. More precisely:
\begin{Theorem} 
\label{smooth}
Let $V$ be a smooth solution to the Euler equations with initial condition $P_0$
and let $A$ be the unique solution of the linear final-value problem
\begin{equation}\label{fvp}
\partial_t A_i+(\delta_{ij}+\partial_j A_i+\partial_i A_j)V^j+\partial_i\varphi=0,
\;\;\;\partial^j A_j=0,\;\;\;A(T,\cdot)=0.
\end{equation}
If  $T$ is small enough the matrix-valued field $\delta_{ij}+\partial_j A_i+\partial_i A_j$
stays uniformly bounded away from zero in $[0,T]\times D$ and, then, the pair $(E,B)$
defined by
$$
E_i=\partial_t A_i+\partial_i\varphi\;,
\;\;\;B_{ij}=
\frac{1}{2}(\partial_jA_i+\partial_i A_j)
$$
is a a maximizer for the concave maximization problem (\ref{ivpdual0}),
while $V$ can be recovered as  $V=-(\mathbb{I}_d+2B)^{-1}E$.
\end{Theorem}
\subsection*{Proof of Theorem \ref{smooth}}

Let $V$ be a smooth solution to the Euler equations, which implies that $V$
is weak solution in the sense of Proposition \ref{weak-new}.
Let us solve the linear final-value problem (\ref{fvp}),
which means that $A$ is a time-dependent divergence-free vector of zero spatial mean,
vanishing at $t=T$ such that, in coordinates,
$$
\partial_t A_i+(\delta_{ij}+\partial_j A_i+\partial_i A_j)V^j+\partial_i \varphi=0,
\;\;{\rm{in}}\;[0,T]\times D,\;\;\int_D A_i(t,\cdot)=0,
\;\;\forall t\in [0,T],
$$
for some scalar field $\varphi$, which can be alternately written
$$
\partial_t A_i+V^j\partial_j A_i-A_j\partial_i V^j+\partial_i \psi+V_i=0,
\;\;\partial^j A_j=0,\;\;A(T,\cdot)=0,\;\;\int_D A_i(t,\cdot)=0,
$$
(where $\psi$ is just $\varphi+A_jV^j$).
Because $V$ is supposed to be smooth, this linear problem can be solved by standard methods
and admits
a unique smooth solution $A$ on $[0,T]\times D$. 
In addition, there is a positive
time $T_0$ depending only on $V$, such that, as long as $0<T\le T_0$,
the field of symmetric matrices 
$\delta_{ij}+\partial_j A_i+\partial_i A_j,$
which is just $\delta_{ij}$
at time $T$,
stays bounded away from $0$,
uniformly in $(t,x)\in [0,T]\times D$.
In this case, we deduce from (\ref{fvp})
\begin{equation}\label{V}
E=-(\mathbb{I}_d+2B)\cdot V,\;\;\;V=-(\mathbb{I}_d+2B)^{-1}\cdot E,
\end{equation}
where $B_{ij}=\partial_j A_i+\partial_i A_j$ and $E_i=\partial_ t A_i+\partial_i\varphi$.
Thus, using notation (\ref{K0}), namely
$$
K[E,B]=\frac{1}{2}\int_{[0,T]\times D}E\cdot(\mathbb{I}_d+2B)^{-1}\cdot E,
$$
we get
$$
2K[E,B]=-\int_{[0,T]\times D}V\cdot E=\int_{[0,T]\times D}V\cdot(\mathbb{I}_d+2B)\cdot V
=\int_{[0,T]\times D}2V\cdot B\cdot V+|V|^2.
$$
By definition, $(E,B)$ are smooth and satisfies constraint (\ref{constraint-B-E}), namely
$$
\partial_t B=LE,\;\;\;B(T,\cdot)=0.
$$
Using that $V$ is a weak solution of the Euler equations in the sense of Proposition (\ref{weak-new}), 
we get (\ref{weakbis}), namely
$$
\int_{[0,T]\times D} V\cdot B\cdot V+(V-P_0)\cdot E=0.
$$
Thus
$$
2K[E,B]=\int_{[0,T]\times D}2P_0\cdot E-2V\cdot E+|V|^2
=4K[E,B]+\int_{[0,T]\times D}2P_0\cdot E+|V|^2,
$$
which shows 
$$
-K[E,B]-\int_{[0,T]\times D}P_0\cdot E=
\frac{1}{2}\int_{[0,T]\times D}|V|^2.
$$
By definition of problems (\ref{ivp1}) and (\ref{ivpdual0}),
we have on one hand
\begin{equation}\label{euler-classic}
\int_{[0,T]\times D} \frac{1}{2}|V|^2\ge \mathcal{I}[P_0]\ge \mathcal{J}[P_0]
\end{equation}
and, on the other hand,
$$
\mathcal{J}[P_0]\ge -K[E,B]-\int_{[0,T]\times D}P_0\cdot E.
$$
So, we conclude that 
$$
\int_{[0,T]\times D} \frac{1}{2}|V|^2=\mathcal{I}[P_0]=\mathcal{J}[P_0]=
-K[E,B]-\int_{[0,T]\times D}P_0\cdot E,
$$
which means that there is no duality gap and that $(E,B)$ is a maximizer of problem (\ref{ivpdual0}),
from which $V$ can be recovered as $V=-(\mathbb{I}_d+2B)^{-1}\cdot E$.
This completes the proof of Theorem \ref{smooth}.

\section{Convex duality and sub-solutions of the Euler equations}
\label{duality0}
Our last result establishes a link between the relaxed 
concave maximization problem (\ref{ivpdual}) and the concept of sub-solution
for the Euler equations, as discussed in \cite{DS} in the context of "Convex Integration".
Here, we limit ourself to sub-solutions that are continuous on $[0,T]\times D$.
\begin{Definition}
\label{def}
We say that a pair of continuous functions $(Z,M)$ 
on $[0,T]\times D$, respectively valued in $\mathbb{R}^d$ and in the space of $d\times d$ symmetric matrices, is a sub-solution to
the Euler equations with initial condition $P_0$ if
\\
i) $M\ge Z\otimes Z$, holds true point-wise in $[0,T]\times D$, in the sense of symmetric matrices;
\\
ii) for every smooth divergence-free vector fields $A$ on $[0,T]\times D$,
of zero spatial mean, that vanishes
at $t=T$ and
each smooth functions $\varphi$ on $[0,T]\times D$, we have
\begin{equation}\label{diff20}
\int_{(t,x)\in[0,T]\times D} \partial_i\varphi(t,x)Z^i(t,x)=0,
\end{equation}
\begin{equation}\label{diff2}
\int_{(t,x)\in[0,T]\times D} \left(\partial_t A_i(t,x)Z^i(t,x)+\partial_i A_j(t,x)M^{ij}(t,x)\right)=-\int_D A_i(0,x)P^i_0(x)dx.
\end{equation}
\end{Definition}
Just as we did for weak solutions, we get:
\begin{Proposition}
\label{sub-new}
$(Z,M)$ is a continuous sub-solution of the Euler equations  with initial condition $P_0$ 
if and only if $M\ge Z\otimes Z$ and
\begin{equation}
\label{subbis}
\int_{[0,T]\times D} M^{ij}B_{ij}+(Z^i-P_0^i)E_i=0,
\end{equation}
for all smooth fields $(E,B)$
respectively valued
in $\mathbb{R}^d$ and in the set of 
symmetric $d\times d$ matrices, that are subject to constraint (\ref{constraint-B-E}),
namely
$$
B(T,\cdot)=0,\;\;\;\partial_t B_{ij}=
\frac{1}{2}(\partial_j E_i+\partial_i E_j)
+\partial_i\partial_j(-\bigtriangleup)^{-1}\partial^k E_k.
$$
\end{Proposition}

Let us now state the main result of this last section:

\begin{Theorem} 
\label{duality}
Let  $P_0$ be a continuous divergence-free vector-field, with zero spatial mean,  on $D$.
Then the infimum of 
$$
\frac{1}{2}\int_{[0,T]\times D} \delta_{ij}M^{ij}
$$
over all continuous sub-solutions $(Z,M)$ of the Euler equations with initial condition $P_0$
(in the sense of Proposition \ref{sub-new}) is finite and just equal to the optimal value of
the relaxed problem (\ref{ivpdual}), namely
$$
\mathcal{J}[P_0]=\max_{(E,B)\in\mathcal{EB}}\;\; -\int_{[0,T]\times D}  P_0\cdot E+
\frac{1}{2}\;E\cdot (\mathbb{I}_d+2B)^{-1}
\cdot E,
$$
where $\mathcal{EB}$ is the class of all Borel measures subject to constraint 
(\ref{constraint-B-E+}).
\end{Theorem}

\subsection*{Proof of Theorem \ref{duality}}

The infimum considered at the beginning of Theorem \ref{duality} can 
be expressed exactly as
$$
\inf_{(Z,M)} K_1(Z,M)+K_2(Z,M)
$$
where $K_1$ and $K_2$ are
the following convex functions  
\begin{equation}\label{K1}
K_1(Z,M)=
\int_{[0,T]\times D} \frac{1}{2}\delta_{ij}M^{ij}
\;\;{\rm{or}}\;\;+\infty
\end{equation}
whether or not $M\ge Z\otimes Z$ is satisfied point-wise, in the sense of symmetric matrices, 
\begin{equation}\label{K2}
K_2(Z,M)=K_3(Z-P_0,M)
\end{equation}
\begin{equation}\label{K3}
K_3(Z,M)=
\sup\;\{\int_{[0,T]\times D} E_iZ^i+B_{ij}M^{ij},\;\;\;(E,B)\in \mathcal{EB}_{smooth}\},
\end{equation}
where 
$\mathcal{EB}_{smooth}$ denotes the class
of all smooth pairs 
$(E,B)$
defined on $[0,T]\times D$ 
with values respectively taken in $\mathbb{R}^d$ and in the set of $d\times d$
symmetric matrices, which satisfy constraint (\ref{constraint-B-E}), namely
$$
\partial_t B=LE,\;\;\;
B(T,\cdot)=0,
$$
where $L$ is the linear operator defined by (\ref{L}). Notice that $B$ is identically trace-free
(since $\delta^{ij}L^k_{ij}=0$ and $B(T,\cdot)=0$) which implies
$
K_3(Z,M)=0,
$
whenever $Z=0$ and $M$ is identically diagonal, a property that we will use in a moment.
\\
\\
There is a rather obvious point $(Z,M)$, namely 
$$
Z(t,x)=P_0(x),\;\;M(t,x)=(1+|P_0(x)|^2)
\mathbb{I}_d,\;\;\forall (t,x)\in [0,T]\times D,
$$
at which $K_1$ is finite and continuous while $K_2$ is just finite.
\\
\\
$[$Indeed, on one hand, 
$\xi\cdot(M-Z\otimes Z)\cdot\xi=(1+|P_0|^2)|\xi|^2-(P_0\cdot \xi)^2\ge |\xi|^2$, 
for all $\xi\in\mathbb{R}^d$,
which implies that $M\ge Z\otimes Z$ holds true in a neighborhood (for the sup-norm) of $(Z,M)$
and, therefore, by definition (\ref{K1}), 
$$
K_1(Z,M)=\int_{[0,T]\times D} \frac{1}{2}\delta_{ij}M^{ij}
$$
is finite and continuous in that neighborhood. 
On the other hand,
$$
K_2(Z,M)=K_3(Z-P_0,M)=K_3(0,M)=0,
$$
since $M$ is diagonal, as just noticed above.$]$
\\
\\
 This is enough to apply Rockafellar's duality theorem (as stated in Chapter 1
of Brezis' book \cite{Brz}) and deduce
\begin{equation}
\label{rocka}
\inf_{(Z,M)} K_1(Z,M)+K_2(Z,M)
=\max_{(E,B)} -K^*_1(-E,-B)-K^*_2(E,B),
\end{equation}
where the minimization is performed over the Banach space of all pairs of continuous functions 
$(Z,M)$ valued in $\mathbb{R}^d$ and in the space of $d\times d$ symmetric matrices, while
the maximization is performed over the dual Banach space of all pairs of bounded Borel
measures $(E,B)$, still valued
in $\mathbb{R}^d$ and in the space of $d\times d$ symmetric matrices.
In this relation,
$K^*_1$, $K^*_2$ denote the Legendre-Fenchel transform of $K_1$ and $K_2$, namely
$$
K^*_1(E,B)=\sup_{(Z,M)}-K_1(Z,M)+\int_{[0,T]\times D}  Z^iE_i+M^{ij}B_{ij}
$$
$$
K^*_2(E,B)=\sup_{(Z,M)}-K_2(Z,M)+\int_{[0,T]\times D}  Z^iE_i+M^{ij}B_{ij}.
$$
In the duality equality (\ref{rocka}), 
the calculation of $K^*_1(-E,-B)$ is easy. Indeed, 
$$
K^*_1(-E,-B)=\sup_{M\ge Z\otimes Z}\int_{[0,T]\times D}  -Z^iE_i-\frac{1}{2}M^{ij}
(\delta_{ij}+2B_{ij})
$$
is equal to $+\infty$, unless $\mathbb{I}_d+2B
\ge 0$ in the sense of symmetric matrices, in which case
$$
K^*_1(-E,-B)=\sup_{Z}\int_{[0,T]\times D}  -Z^iE_i-\frac{1}{2}(\delta_{ij}+2B_{ij})Z^i Z^j,
$$
which we can write in the more compact (but less precise) form:
$$
K^*_1(-E,-B)=\frac{1}{2}\int_{[0,T]\times D}  E\cdot (\mathbb{I}_d+2B)^{-1}\cdot E\;\;\in [0,+\infty].
$$
Let us now compute $K_2^*$. We find (by definition (\ref{K2})
$$
K^*_2(E,B)=\sup_{(Z,M)}\;-K_3(Z-P_0,M)+\int_{[0,T]\times D}  E_i Z^i+B_{ij}M^{ij}
$$
$$
=K_3^*(E,B)+\int_{[0,T]\times D}  E_i P_0^i.
$$
By definition of $K_3$ in (\ref{K3}), $K^*_3(E,B)$ is just zero or $+\infty$ whether or not
$(E,B)$ belongs to the weak-* closure 
of class $\mathcal{EB}_{smooth}$.
This closure is nothing but the class $\mathcal{EB}$
of all bounded Borel measures $(E,B)$ valued
in $\mathbb{R}^d$ and in the space of $d\times d$ symmetric matrices that satisfy constraint
(\ref{constraint-B-E}), or, more precisely (\ref{constraint-B-E+}). 
$[$The reason being the same as the one for which $\mathcal{EB}_{\infty,2}$ is the 
weak-* closure of class $\mathcal{EB}_{smooth}$ in $L^2\times L^\infty$ as already explained
when proving Proposition \ref{weak-new} in section \ref{reformulation}.$]$
Therefore, the duality equality (\ref{rocka}) exactly gives the result stated in Theorem 
\ref{duality} which completes the proof.


\end{document}